\documentclass[12pt]{article} 
\usepackage{amsmath,amsxtra,amssymb,latexsym, amscd,amsthm}

\advance\voffset-1truecm\relax
\advance\hoffset-1.0truecm\relax
\font\titbf=cmbx10 scaled \magstep2

\font\tac=cmcsc10 scaled \magstep1

\numberwithin{equation}{section}

\newtheorem{theorem}{Theorem}
\newtheorem{corollary}[theorem]{Corollary}
\newtheorem{lemma}[theorem]{Lemma}
\def\C{\Bbb C}
\def\Z{\Bbb Z}

\begin{document}
$\qquad$

\centerline{\titbf  UNIQUENESS THEOREMS FOR MEROMORPHIC }
\medskip

\centerline{\titbf  MAPPINGS WITH FEW HYPERPLANES}
\medskip

\centerline{\tac Gerd Dethloff and Tran Van Tan }

\vskip0.15cm

\begin{abstract}
\noindent Let $f, g$ be linearly nondegenerate meromorphic mappings  of $\mathbb{C}^m$ into $\mathbb{C} P^n.$  Let $\{H_j\}_{j=1}^q$  be hyperplanes in $\mathbb{C} P^n$ in general position, such that 

$\text{a)}\quad$   $f^{-1}(H_j)= g^{-1}(H_j)\ ,\quad \text{for all}\ \ 1\leq j\leq  q,$

$\text{b)}\quad$   $\dim \big(f^{-1}(H_i)\cap f^{-1}(H_j)\big)\le m-2\ \ \text{for all}\  1\le i<j\le q,$ and

$\text{c)}\quad$  $f=g$ on \ \
$\bigcup_{j=1}^{q}f^{-1}(H_j).$

\noindent It is well known that if $q\geq 3n+2$, 
then $f\equiv g.$ In this paper we show that for every nonnegative integer $c$ there exists positive integer $N(c)$ depending only on   $c$ in an explicit way  such that the above result remains valid if $q\geq (3n+2 -c)$ and $n\geq N(c).$ Furthermore, we also show that the coefficient of $n$ in the formula of $q$ can be replaced by a number which is strictly smaller than 3 for all $n>>0.$ At the same time, a big number of recent uniqueness theorems are generalized considerably.


\end{abstract}
\section{Introduction}
The uniqueness problem of meromorphic mappings under a condition
on the inverse images of divisors was first studied by R. Nevanlinna \cite{Ne}. He
showed that for two nonconstant meromorphic functions $f$ and $g$ on the
complex plane $\mathbb{C}$, if they have the same inverse images for five
distinct values, then $f\equiv g.$ We remark that the number of distinct values in the above result can not be replaced by a smaller one, as it can be seen easily as follows: Let $f$ be a nonconstant nonvanishing holomorphic function on $\C,$ then consider the two distinct functions $f, \frac{1}{f}$ and the four values $0, \infty, 1,-1.$ In 1975, H. Fujimoto \cite{F1} generalized
R. Nevanlinna's result to the case of meromorphic mappings of $\mathbb{C}^{m}$
into $\mathbb{C}P^{n}$. He showed  that for two linearly nondegenerate meromorphic
mappings $f$ and $g$ of $\mathbb{C}^{m}$ into $\mathbb{C}P^{n}$, if they
have the same inverse images counted with multiplicities for $(3n+2)$
hyperplanes in general position in $\mathbb{C}P^{n},$ then $f\equiv g$. Since that time, this problem has been  studied intensively by H. Fujimoto, W. Stoll, L. Smiley, S. Ji, M. Ru, G. Dethloff, T. V. Tan, D. D. Thai, S. D. Quang and others.

In 1983, L. Smiley \cite{Sm} showed that
\begin{theorem}
 Let $f, g$ be linearly nondegenerate meromorphic mappings  of $\mathbb{C}^m$ into $\mathbb{C} P^n.$  Let $\{H_j\}_{j=1}^q$  ($q\geq 3n+2)$ be hyperplanes in $\mathbb{C} P^n$ in general position. Assume that

$\text{a)}\quad$   $f^{-1}(H_j)= g^{-1}(H_j)\ ,\quad \text{for all}\ \ 1\leq j\leq  q,$ (as sets)

$\text{b)}\quad$   $\dim \big(f^{-1}(H_i)\cap f^{-1}(H_j)\big)\le m-2\ \ \text{for all}\  1\le i<j\le q,$ and

$\text{c)}\quad$  $f=g$ on \ \
$\bigcup_{j=1}^{q}f^{-1}(H_j).$

\noindent Then $f\equiv g.$
\end{theorem}
\noindent Theorem 1 was given again in 1989 by W. Stoll \cite{St} and in 1998 by H. Fujimoto \cite{F2}.   There is a number of papers which tried to extend Theorem 1 to the case of fewer hyperplanes. 
For example, in 1988 S. Ji \cite{Ji} considered three  linearly nondegenerate meromorphic mappings  $f, g, h$ of $\mathbb{C}^m$ into $\mathbb{C} P^n,$ and he showed that if for any two mappings of them the conditions  $a), b), c)$ are satisfied, then $f\times g\times h$ is algebraically degenerate. In 2006 in \cite{DT3} we showed that the result of S. Ji remains valid in the case $q\geq [\frac{5(n+1)}{2}]$
 (where we denote $[x]:= \max\{ k\in \Z: k\leq x\}$ for a real number $x$).

 In \cite{DT2} we showed that Theorem 1  remains valid for $n\geq 2$ and $q\geq 3n+1$ hyperplanes, but  the condition $a)$ is replaced by 
$$a'): \min\{\nu_{(f,H_j)},1\}=\min\{\nu_{(g,H_j)},1\}\ (1\leq j\leq  2n-2),\text{ and\ } $$
$$\min\{\nu_{(f,H_j)},2\}=\min\{\nu_{(g,H_j)},2\} \ (2n-1\leq j\leq  3n+1).$$ 
But in all of these results either the assertion is weaker (i.e. one did not get $f \equiv g$) or the assumption is stronger,
in the sense that the conditions $f^{-1}(H_j)= g^{-1}(H_j)$ do not only hold set-theoretically, but with counting multiplicities, at least up to a certain order (we refer the reader to \cite{DT2} - \cite{F2}, \cite{Ji}, \cite{St}, \cite{TQ2} for further results and comments on it).  The only exception seems to be the recent result of Thai and Quang \cite{TQ2}, which slightly improves our result in 
\cite{DT2} mentioned above, by proving it only under the original condition $a)$ instead of
the condition $a')$, and, thus, gives a generalization of Smiley's Theorem 1 in the strict sense:
\begin{theorem}
Let $n \geq 2$ and $f, g$, $\{H_j\}_{j=1}^q$ be as in  Theorem 1. Then for $q\geq 3n+1$, one has $f\equiv g.$
\end{theorem}
\noindent In the same paper, they asked if for $q < 3n+1$, there exist positive
integers $N_0$ such that for $n \geq N_0$, the above unicity theorems hold.

In this paper we show that for every nonnegative integer $c$ there exists a positive integer $N(c)$ depending only on   $c$ such that the above unicity theorems Theorem 1 and Theorem 2 remain valid if $q\geq (3n+2 -c)$ and $n\geq N(c)$. We also get that the coefficient of $n$ in the formula of $q$ can be replaced by a number which is smaller than 3 for all $n>>0.$ Thus, we get affirmative answers to
the question of Thai and Quang.
But our main result Theorem 3 below  is in fact much stronger,
it does not only improve considerably Theorem 1 and Theorem 2, but also many other uniqueness theorems, taking into account (truncated) orders of the inverse images of the hyperplanes:

\begin{theorem} \label{E:pt3}
Let $f$ and $g$ be two linearly nondegenerate meromorphic mappings of $\C^m$ into $\C P^n$ and let
$H_1, \dots, H_q$ $(q \geqslant 2n)$ be hyperplanes in $\C P^n$
in general position. Let $p$ be a positive integer. Assume that

$\text{a)}\quad$ $\min\{\nu_{(f,H_j)},p\}=\min\{\nu_{(g,H_j)},p\}$ for all  $1\le j\le q,$

$\text{b)}\quad$ $\dim \big(f^{-1}(H_i)\cap f^{-1}(H_j)\big)\le m-2$ for all $  1\le i<j\le q,$ and

$\text{c)}\quad$  $f=g$ on $\bigcup_{j=1}^{q}f^{-1}(H_j).$

\noindent Then the following assertions hold:

$1)\quad$ If $p=n, q\geq 2n+3$ then $f\equiv g.$

$2)\quad$ If $p<n$ and there exists a positive integer $t\in\{p,\cdots,n-1\}$ such that
$$\big(\frac{(q+2t)(q-n-1)}{n}-2q\big)\frac{(n-t)(q+2p-2)}{4nt}>2q-\frac{(q-n-1)(q+2p-2)}{n}$$
 then $f\equiv g.$
\end{theorem}

\noindent {\bf Remark.}
a) The assertion 1) of Theorem 3 is a kind of generalization of R. Nevannina's result to the case of meromorphic mappings of $\C^m$ into $\C P^n$. \\
b) Theorem 3 gives also the solution for the open questions which were given by H. Fujimoto in \cite{F2}, \cite{F3}.

The most interesting special cases of Theorem 3 are  the cases $p=n$ and  $p=1$.

The case $p=n$ is the one which gives the unicity theorem with the fewest number of hyperplanes, namely $2n+3$. The best result known before was
our result in \cite{DT2}, where we showed that, under the same assumptions, the unicity theorem holds for 
$n \geq 2$ and $q \geq n + [ \sqrt{2n(n+1)}]+2$.

The case $p=1$ is the one  where multiplicities of the inverse images of the hyperplanes are not taken into account as in the Theorems 1 and 2 of Smiley and Thai-Quang. In this case, the inequality in the assertion 2) of Theorem 3 will become the following 
$$(*)\quad\big(\frac{(q+2t)(q-n-1)}{n}-2q\big)\frac{(n-t)}{4nt}>2-\frac{q-n-1}{n}.$$ 
We state some cases where the condition $(*)$ is satisfied
(we remark that in all these cases the condition $1 \leq t \leq n-1$ is satisfied):

1)\quad $n\geq 2, q\geq 3n+1, t=1.$

2)\quad $n\geq 5, q\geq 3n, t=2.$

3)\quad $n\geq 5(c+1)^2, q=3n-c, t=3c$  for each positive integer $c \geq 1.$ 

4)\quad If $t=[\frac{n}{2}],q=[\frac{11n}{4}] $, then the right side of $(*)$ $\geq\frac{17}{64}-O(\frac{1}{n})$, and the left side of $(*)$ $\leq \frac{1}{4}+O(\frac{1}{n}).$ So, in this case $(*)$ is satisfied for all $n>>0.$ \\
Case 1) gives Theorem 2 of Thai-Quang above, and cases 2) and 3) lead  to the following:

\begin{corollary}
 Let $c \geq 0$ a given non negative integer. Let $n \geq 5(c+1)^2$
 and $q=3n-c$. 

 Let $f, g$ be linearly nondegenerate meromorphic mappings  of  $\mathbb{C}^m$ into $\mathbb{C} P^n.$  Let $\{H_j\}_{j=1}^q$   be hyperplanes in $\mathbb{C} P^n$ in general position. Assume that

$\text{a)}\quad$   $f^{-1}(H_j)= g^{-1}(H_j)\ ,\quad \text{for all}\ \ 1\leq j\leq  q,$ (as sets)

$\text{b)}\quad$   $\dim \big(f^{-1}(H_i)\cap f^{-1}(H_j)\big)\le m-2\ \ \text{for all}\  1\le i<j\le q,$ and

$\text{c)}\quad$  $f=g$ on \ \
$\bigcup_{j=1}^{q}f^{-1}(H_j).$\\
\noindent Then $f\equiv g.$

\end{corollary}

\noindent Case 4) leads to the following:
\begin{corollary}
 There exists a natural number $n_0>>0$ such that for  $ n \geq n_0$ and $q=[2,75n]$
 the following holds:

 Let $f, g$ be linearly nondegenerate meromorphic mappings  of $\mathbb{C}^m$ into $\mathbb{C} P^n.$  Let $\{H_j\}_{j=1}^q$   be hyperplanes in $\mathbb{C} P^n$ in general position. Assume that

$\text{a)}\quad$   $f^{-1}(H_j)= g^{-1}(H_j)\ ,\quad \text{for all}\ \ 1\leq j\leq  q,$ (as sets)

$\text{b)}\quad$   $\dim \big(f^{-1}(H_i)\cap f^{-1}(H_j)\big)\le m-2\ \ \text{for all}\  1\le i<j\le q,$ and

$\text{c)}\quad$  $f=g$ on \ \
$\bigcup_{j=1}^{q}f^{-1}(H_j).$\\
\noindent Then $f\equiv g.$

\end{corollary}

\noindent {\bf Remark.}
The number $n_0$ can be explicitly calculated. 

Let us finally give some comments on our method of proof. First of all, we use  an auxiliary function for estimating the counting function, which is different from the auxiliary functions which were used in the previous papers. Thanks to this point, the estimate which we obtain here is better than the estimate of the previous authors (including ourselves) if $p>1.$ After that, we try to replace the value at which multiplicities are truncated by a bigger one. This idea did not appear in the previous papers. In order to carry out this idea, we estimate the counting function of the set $A$ of all points with multiplicites in $\{p,\cdots,t\}.$ Then combining with the assumption "multiplicites are truncated by $p$" , we see the condition  "multiplicites are truncated by $t+1$" is satisfied automatically outside $A.$ Thanks to this technique, if $p<n,$ we will get a stronger version for the Second Main Theorem for meromorphic mappings $f$ and $g$ with hyperplanes $\{H_j\}_{j=1}^q.$ Hence, with this method we will get better uniqueness theorems if $p>1$ or $p<n.$ This means that we get a better uniqueness theorem unless  $n=p=1.$ This perfectly coincides with the fact that the result of R. Nevanlinna is optimal as we remarked above.

\vskip0.5cm
\noindent
{\bf Acknowledgements:} The second author would like to thank the Universit\'{e} de Bretagne Occidentale for its hospitality and for support, PICS-CNRS ForMath Vietnam for support.

\section{Preliminaries}

We set $\Vert z \Vert := \big(|z_1|^2 + \cdots + |z_m|^2\big)^{1/2}$
for $z = (z_1, \dots, z_m) \in \C^m$ and define
$$B(r) := \big\{z \in \C^m : \Vert z \Vert < r\big\},\qquad
S(r) := \big\{z \in \C^m : \Vert z \Vert = r\big\}
$$ 
for all $0 < r < \infty$.
Define
\begin{align*}
d^c &:= \frac{\sqrt{-1}}{4\pi}(\overline \partial - \partial),\quad
\upsilon := \big(dd^c\Vert z \Vert^2 \big)^{m-1}\quad \\
\sigma &:= d^c \text{log}\Vert z\Vert^2 \land 
\big(dd^c \text{log}\Vert z\Vert^2\big)^{m-1}.
\end{align*}
Let $F$ be a nonzero holomorphic function on $\C^m$.
For each $a \in \C^m$, expanding $F$ as $F = \sum P_i(z-a)$
with homogeneous polynomials $P_i$ of degree $i$ around $a$,
we define
$$ \nu_F(a) := \min \big\{ i : P_i \not\equiv 0 \big\}.$$
Let  $\varphi$ be a nonzero meromorphic function on $\C^m$.
We define the zero divisor $\nu_\varphi$ as follows: For each $z \in 
\C^m$, we choose nonzero holomorphic functions $F$ and $G$ on a 
neighborhood $U$ of $z$ such that $\varphi = {F}/{G}$ on
$U$ and $\text{dim}\,\big(F^{-1}(0) \cap G^{-1}(0)\big) \leqslant m-2$.
Then we put $\nu_\varphi(z) := \nu_F(z)$.

\noindent Let $\nu$ be a divisor in $\mathbb{C}^m$ and $k$ be positive integer or $+\infty $. Set \ $
|\nu|:=\overline{\big\{z:\ \nu(z)\neq 0
\big\}}$  and $\nu^{[k]}(z) := \min \{ \nu_\varphi(z), k\}.$ 

The truncated counting function of $\nu$ is defined by
$$ N^{[k]}(r, \nu) := \int\limits_1^r 
\frac{n^{[k]}(t)}{t^{2m-1}} dt \quad (1 < r < + \infty), $$
where
\begin{align*}
n^{[k]}(t) = \begin{cases}
\displaystyle{\int\limits_{|\nu| \cap B(t)}}
\nu^{[k]}\cdot \upsilon \ &\text{for}\ m \geqslant 2,\\
\sum\limits_{|z| \leqslant t} \nu^{[k]}(z)
&\text{for}\ m = 1.\end{cases}
\end{align*}
We simply write \ $N(r, \nu)$ for $N^{[+\infty]}(r,\nu)$.

\noindent For a nonzero meromorphic function $\varphi$ on $\mathbb{C}^{m},$ we set \quad $N_{\varphi}^{[k]}(r):=N^{[k]}(r, \nu_\varphi)$ and $ N_{\varphi}(r):=N^{[+\infty]}(r, \nu_\varphi).$ 
We have the following Jensen's formula:
$$ N_\varphi (r) - N_{\frac{1}{\varphi}}(r) =
\int\limits_{S(r)} \text{log}|\varphi| \sigma 
- \int\limits_{S(1)} \text{log}|\varphi| \sigma .$$

For a closed subset $A$ of a purely $(m-1)$-dimensional analytic subset 
of $
\mathbb{C}^{m}$, we define
\begin{equation*}
N^{[1]}(r,A):=\int\limits_{1}^{r}\frac{n^{[1]}(t)}{t^{2m-1}
}dt,\;\;\;(1<r<+\infty ) 
\end{equation*}
where
\begin{equation*}
n^{[1]}(t):=\begin{cases}\int\limits_{A\cap B(t)}\upsilon 
\;\;\;\;\;\;\;\;\; \text{  for }m\geq 2\\
\sharp\left( A\cap B(t)\right) \text{ for }m=1.\end{cases} 
\end{equation*}

Let $f : \C^m \longrightarrow \C P^n$ be a meromorphic mapping.
For an arbitrary fixed homogeneous coordinate system 
 $(w_0 : \cdots : w_n)$ in
$\C P^n$, we take a reduced representation $f = (f_0 : \cdots : f_n)$,
which means that each $f_i$ is a holomorphic function on $\C^m$
and $f(z) = (f_0(z) : \cdots : f_n(z))$ outside the analytic set
$\{ f_0 = \cdots = f_n = 0\}$ of codimension $\geqslant 2$.
Set $\Vert f \Vert = \big(|f_0|^2 + \cdots + |f_n|^2\big)^{1/2}$.
The characteristic function $T_f(r)$ of $f$ is defined by
$$ T_f(r) := \int\limits_{S(r)} \text{log} \Vert f \Vert \sigma -
\int\limits_{S(1)} \text{log} \Vert f \Vert \sigma , \quad
1 < r < + \infty . $$   
For a meromorphic function $\varphi$ on $\C^m$, the characteristic
function $T_\varphi(r)$ of $\varphi$ is defined by considering
$\varphi$ as a meromorphic mapping of $\C^m$ into $\C P^1$.

\noindent The proximity function $m(r,\varphi)$ is defined by
$$ m(r,\varphi) = \int\limits_{S(r)} \text{log}^+ |\varphi| \sigma, $$
 where $\text{log}^+ x = \max \big\{ \text{log}x, 0\big\}$ for 
$x \geqslant 0$.

We state the First and Second Main Theorem in Value Distribution 
Theory:

\noindent
For a hyperplane $H : a_0 w_0 + \cdots + a_n w_n = 0$ in $\C P^n$
with $f(\C^m) \not\subseteq H$, we put $(f,H) = a_0 f_0 + \cdots 
+ a_n f_n$, where $(f_0 : \cdots : f_n)$  is a reduced
representation of $f$.\\
\medskip
\noindent
{\bf First Main Theorem.} 

1) {\it For a nonzero meromorphic function $\varphi,$ on $\C^m$ we have }
$$ T_\varphi (r) = N_{\frac{1}{\varphi}}(r) + m(r,\varphi) + O(1).$$

2) {\it Let $f$ be a meromorphic mapping of
$\C^m$ into $\C P^n$, and $H$ be a hyperplane in $\C P^n$ such that $(f,H) \not\equiv 0$. Then}
$$ N_{(f,H)}(r) \leqslant T_f(r) + O(1) \quad \text{\it for all}\ r > 
1.$$

\medskip
\noindent
{\bf Second Main Theorem.} {\it Let $f$ be a linearly nondegenerate
meromorphic mapping of $\C^m$ into $\C P^n$ and
$H_1, \dots, H_q$ $(q \geqslant n + 1)$ hyperplanes in $\C P^n$
in general position. Then
$$ (q-n-1) T_f(r) \leqslant \sum_{j=1}^q N_{(f,H_j)}^{[n]}(r) +
o\big(T_f(r)\big) $$ 
for all $r$ except for a subset $E$ of $(1, +\infty)$ of finite 
Lebesgue measure.}

\section{Proof of Theorem 3 }
In order to prove Theorem 3 we need the following lemma.
\begin{lemma}\label{Le3.1}Let $f$ and $g$ be two distinct linearly nondegenerate mappings of $\C^m$ into  $\C P^n$ 
and let $H_1, \dots, H_q$ $(q \geqslant n+1)$ be hyperplanes in $\C P^n$
in general position. Assume that

$\text{a)}\quad$ $\min\{\nu_{(f,H_j)},1\}=\min\{\nu_{(g,H_j)},1\}$ for all  $1\le j\le q,$

$\text{b)}\quad$ $\dim \big(f^{-1}(H_i)\cap f^{-1}(H_j)\big)\le m-2$ for all $  1\le i<j\le q,$ and

$\text{c)}\quad$  $f=g$ on $\bigcup_{j=1}^{q}f^{-1}(H_j).$

 \noindent Then for every positive integer $\ell$ and for every subset $\{i_0, j_0\}\subset\{1,\cdots,q\}$ such that
 \begin{equation*}
\text{det} 
\begin{pmatrix}
(f,H_{i_0}) &  & (f,H_{j_0})\cr&  \cr (g,H_{i_0}) &  & 
(g,H_{j_0})
\end{pmatrix}
\not\equiv 0
, \text{\ \ we have}
\end{equation*}
\begin{eqnarray*}
\sum_{j=1,j\ne i_{0},j_{0}}^{q}N
_{(f,H_j)}^{[1]}(r) &+N_{(f,H_{i_0})}^{[\ell]}(r)+ N_{(f,H_{j_0})}^{[\ell]}(r) \leq 
T_f(r)+T_g(r)\\
&+(\ell-1)\big(N^{[1]}(r,\overline{A})
+N^{[1]}(r,\overline{B})\big)+O(1),
\end{eqnarray*}
where $A:=\{z: \min\{\nu_{(f,H_{i_0})}(z),\ell\}\ne\min\{\nu_{(g,H_{i_0})}(z),\ell\}\}$, and
$B:=\{z: \min\{\nu_{(f,H_{j_0})}(z),\ell\}\ne\min\{\nu_{(g,H_{j_0})}(z),\ell\}\}.$
\end{lemma}
\noindent
\textbf{Proof.}
Set $$\phi:=\frac{(f,H_{i_0})}{(f,H_{j_0})}-\frac{(g,H_{i_0})}{(g,H_{j_0})}\not\equiv 0.$$
Let $z_0$ be an arbitrary zero point of $(f,H_{i_0})$ (if there exist any). If $z_0\in A$, then $z_0$ is a zero point of  $\phi$ (outside an analytic set of codimension $\geq 2$). If $z_0\not\in A,$ then we have 
$$\min\{\nu_{(f,H_{i_0})}(z_0),\ell\}=\min\{\nu_{(g,H_{i_0})}(z_0),\ell\}. $$
\noindent In this case, $z_0$ is a zero point of $\phi$ with multiplicity $\geq \min\{\nu_{(f,H_{i_0})}(z_0),\ell\}$ (outside an analytic set of codimension $\geq 2$). 

\noindent For any $j\in\{1,\cdots,q\}\backslash\{i_0,j_0\},$ since $f=g$ on $f^{-1}(H_j)$ we have that a zero point of $(f,H_j)$ is also a zero point of $\phi$ (outside an analytic set of codimension $\geq 2$). 

\noindent On the other hand $\dim \big(f^{-1}(H_i)\cap f^{-1}(H_j)\big)\le m-2$ for all $  1\le i<j\le q.$ 
\noindent Hence, we have

\begin{align}
N_\phi(r)\geq N_{(f,H_{i_0})}^{[\ell]}-(\ell-1) N^{[1]}(r,\overline{A})+ \sum_{j=1,j\ne i_{0},j_{0}}^{q}N
_{(f,H_j)}^{[1]}(r). \label{1}
\end{align}
By the First Main Theorem we have
\begin{align*}
m\big(r,\frac{(f,H_{i_0})}{(f,H_{j_0})}\big)&=T_{\frac{(f,H_{i_0})}{(f,H_{j_0})}}(r)-N_{\frac{(f,H_{j_0})}{(f,H_{i_0})}}(r)+O(1)\\
&\leq T_f(r)-N_{(f,H_{j_0})}(r) +O(1).
\end{align*}
Similarly (note that $f^{-1}(H_j) = g^{-1}(H_j)\, , \: j=1, \dots , q$ by condition $a)$ of Lemma 6),
\begin{align*}
m\big(r,\frac{(g,H_{i_0})}{(g,H_{j_0})}\big)\leq T_g(r)-N_{(g,H_{j_0})}(r) +O(1).
\end{align*}
Hence, we have
\begin{align}
m(r,\phi)&\leq m\big(r,\frac{(f,H_{i_0})}{(f,H_{j_0})}\big)+m\big(r,\frac{(g,H_{i_0})}{(g,H_{j_0})}\big)+O(1)\notag\\
&\leq T_f(r)+T_g(r)-N_{(f,H_{j_0})}(r)-N_{(g,H_{j_0})}(r) +O(1).\label{2}
\end{align}

Set $\nu=\max\{\nu_{(f,H_{j_0})},\nu_{(g,H_{j_0})}\}.$

\noindent It is clear that 
\begin{align}
\nu+\nu_{(f,H_{j_0})}^{[\ell]}- \nu_{(f,H_{j_0})}-\nu_{(g,H_{j_0})}\leq \ell -1 \text{\ on \ }  B\label{3}
\end{align}
(note that $f^{-1}(H_{j_0})=g^{-1}(H_{j_0})).$

\noindent Since $\min\{\nu_{(f,H_{j_0})},\ell\}=\min\{\nu_{(g,H_{j_0})},\ell\}$ on $\C^m\backslash B$ we have
\begin{align}
\nu+\nu_{(f,H_{j_0})}^{[\ell]}- \nu_{(f,H_{j_0})}-\nu_{(g,H_{j_0})}\leq 0 \text{\ on \ } \C^m\backslash B\label{4}.
\end{align}
By (\ref{3}), (\ref{4}) we have
\begin{align*}
N_{(f,H_{j_0})}(r)+N_{(g,H_{j_0})}(r)\geq N(r,\nu)+ N_{(f,H_{j_0})}^{[\ell]}(r)-(\ell-1)N^{[1]}(r,\overline B).
\end{align*}
Combining with (\ref{2}) we have
\begin{align*}
m(r,\phi)\leq T_f(r)+T_g(r)-N(r,\nu)- N_{(f,H_{j_0})}^{[\ell]}(r)+(\ell-1)N^{[1]}(r,\overline B)+O(1).
\end{align*}
On the other hand, it is clear that
\begin{align*}
N(r,\nu)\geq N_{\frac{1}{\phi}}(r).
\end{align*}
Hence, we get
\begin{align*}
m(r,\phi)\leq T_f(r)+T_g(r)-N_{\frac{1}{\phi}}(r)- N_{(f,H_{j_0})}^{[\ell]}(r)+(\ell-1)N^{[1]}(r,\overline B)+O(1).
\end{align*}
Then, by the First Main Theorem we have
\begin{align}
N_{\phi}(r)&\leq  T_\phi(r)+O(1)=m(r,\phi)+N_{\frac{1}{\phi}}(r)+O(1)\notag\\
&\leq T_f(r)+T_g(r)- N_{(f,H_{j_0})}^{[\ell]}(r)+(\ell-1)N^{[1]}(r,\overline B)+O(1).\label{5}
\end{align}
By (\ref{1}) and (\ref{5}) we have
\begin{align*}
N_{(f,H_{i_0})}^{[\ell]}&-(\ell-1) N^{[1]}(r,\overline{A})+ \sum_{j=1,j\ne i_{0},j_{0}}^{q}N
_{(f,H_j)}^{[1]}(r)\\
&\leq T_f(r)+T_g(r)- N_{(f,H_{j_0})}^{[\ell]}(r)+(\ell-1)N^{[1]}(r,\overline B)+O(1).
\end{align*}
This gives 
\begin{align*}
\sum_{j=1,j\ne i_{0},j_{0}}^{q}&N
_{(f,H_j)}^{[1]}(r)+N_{(f,H_{i_0})}^{[\ell]}+N_{(f,H_{j_0})}^{[\ell]}(r)\\
&\leq T_f(r)+T_g(r)+(\ell-1) N^{[1]}(r,\overline{A})+(\ell-1)N^{[1]}(r,\overline B)+O(1).
\end{align*}

We have completed proof of Lemma 6.
\hfill$\Box$
\vskip0.5cm
\noindent\textbf{Proof of Theorem 3.}
Assume that $f\not\equiv g.$

\noindent We introduce an equivalence relation on $ L:=\{1,\cdots, q\}$ as follows: $i\sim j$ if and only if
\begin{equation*}
\text{det} 
\begin{pmatrix}
(f,H_{i}) &  & (f,H_{j})\cr&  \cr (g,H_{i}) &  & 
(g,H_{j})
\end{pmatrix}
\equiv 0.
\end{equation*}
Set $\{L_1,\cdots, L_s\}=L/\sim $. Since $f\not\equiv g$ and $\{H_j\}_{j=1}^q$ are in general position, we have that 
 $\sharp L_k\leq n$ for all $k\in\{1,\cdots, s\}.$ Without loss of generality, we may assume that $L_k:=\{i_{k-1}+1,\cdots, i_k\}$
($k\in\{1,\cdots, s\})$ where $0=i_0<\cdots <i_s=q.$ 

\noindent We define the map $\sigma: \{1,\cdots, q\}\to \{1,\cdots, q\}$ by
\begin{equation*}
\sigma (i)=
\begin{cases}
i+n& \text{ if $i+n\leq q$},\\
i+n-q& \text{ if  $i+n>q$}.
\end{cases}
\end{equation*}
It is easy to see that  $\sigma$ is bijective and $\mid \sigma (i)-i\mid\geq n $ (note  that $q\geq 2n).$ This implies that $i$ and $\sigma (i)$ belong two distinct sets of $\{L_1,\cdots, L_s\}.$ This implies that 
\begin{equation*}
\text{det} 
\begin{pmatrix}
(f,H_{i}) &  & (f,H_{\sigma(i)})\cr&  \cr (g,H_{i}) &  & 
(g,H_{\sigma(i)})
\end{pmatrix}
\not\equiv 0.
\end{equation*}
For each $i\in\{1,\cdots, q\},$ by Lemma 6 (with $\ell=p, i_0=i, j_0=\sigma(i))$ we have
\begin{align*}
\sum_{j=1,j\ne i,\sigma(i)}^{q}N
_{(f,H_j)}^{[1]}(r) +N_{(f,H_i)}^{[p]}(r)+ N_{(f,H_{\sigma(i)})}^{[p]}(r) \leq 
T_f(r)+T_g(r)+O(1)
\end{align*}
(note that $\min\{\nu_{(f,H_j)},p\}=\min\{\nu_{(g,H_j)},p\}$ for all  $1\le j\le q).$

\noindent This implies that 
\begin{align*}
(q-2)\sum_{j=1}^{q}N
_{(f,H_j)}^{[1]}(r) +\sum_{i=1}^{q}\big(N_{(f,H_i)}^{[p]}(r)+ N_{(f,H_{\sigma(i)})}^{[p]}(r)\big) \leq 
q\big(T_f(r)+T_g(r)\big)+O(1).
\end{align*}
This gives
\begin{align*}
(q-2)\sum_{j=1}^{q}N
_{(f,H_j)}^{[1]}(r) +2\sum_{i=1}^{q}N_{(f,H_i)}^{[p]}(r) \leq 
q\big(T_f(r)+T_g(r)\big)+O(1).
\end{align*}
Similarly,
\begin{align*}
(q-2)\sum_{j=1}^{q}N
_{(g,H_j)}^{[1]}(r) +2\sum_{i=1}^{q}N_{(g,H_i)}^{[p]}(r) \leq 
q\big(T_f(r)+T_g(r)\big)+O(1).
\end{align*}
Therefore, we get
\begin{align}
(q-2)\sum_{j=1}^{q}\big(N
_{(f,H_j)}^{[1]}(r)&+N
_{(g,H_j)}^{[1]}(r)\big) +2\sum_{i=1}^{q}\big(N_{(f,H_i)}^{[p]}(r)+N_{(g,H_i)}^{[p]}(r)\big) \notag\\
&\leq 2q\big(T_f(r)+T_g(r)\big)+O(1).\label{new1}
\end{align}
By the Second Main Theorem, we have 
(for all $r$ except for a subset $E$ of $(1, +\infty)$ of finite 
Lebesgue measure, which, for simplicity, we do not  mention any more
in the following if no confusion can arise)
\begin{align*}
(q-n-1)(T_f(r)+T_g(r))\leq \sum_{j=1}^{q}\big(N
_{(f,H_j)}^{[n]}(r)+N
_{(g,H_j)}^{[n]}(r)\big)+o(T_f(r)+T_g(r)).
\end{align*}
Hence, by (\ref{new1}) we get
\begin{align}
(q-2)&\sum_{j=1}^{q}\big(N_{(f,H_j)}^{[1]}(r)- \frac{1}{n}N_{(f,H_j)}^{[n]}(r)\big)+ (q-2)\sum_{j=1}^{q}\big
(N_{(g,H_j)}^{[1]}(r)- \frac{1}{n}N_{(g,H_j)}^{[n]}(r)\big)\notag\\
&+2\sum_{j=1}^{q}\big(N_{(f,H_j)}^{[p]}(r)- \frac{p}{n}N_{(f,H_j)}^{[n]}(r)\big)+ 2\sum_{j=1}^{q}\big
(N_{(g,H_j)}^{[p]}(r)- \frac{p}{n}N_{(g,H_j)}^{[n]}(r)\big)\notag\\
&\leq (q-2)\sum_{j=1}^{q}\big(N_{(f,H_j)}^{[1]}(r)+N_{(g,H_j)}^{[1]}(r)\big)+2\sum_{j=1}^{q}\big(N_{(f,H_j)}^{[p]}(r)+N_{(g,H_j)}^{[p]}(r)\big)\notag\\
&\quad-\frac{(q-n-1)(q+2p-2)}{n}(T_f(r)+T_g(r))+o(T_f(r)+T_g(r))\notag\\
&\overset{(\ref{new1})}{\leq}\big(2q-\frac{(q-n-1)(q+2p-2)}{n}\big)(T_f(r)+T_g(r))+o(T_f(r)+T_g(r))\label{8}.
\end{align}
\indent\hspace*{0.3cm}1) If $p=n$ and $q\geq2n+3,$ then $2q-\frac{(q-n-1)(q+2p-2)}{n}<0.$ This contradicts to (\ref{8}). So, we have $f\equiv g.$\hfill$\Box$

\indent\hspace*{0.3cm}2) Assume that $p<n$ and that there exists a positive integer $t\in\{p,\cdots,n-1\}$ such that
\begin{align}
\big(\frac{(q+2t)(q-n-1)}{n}-2q\big)\frac{(n-t)(q+2p-2)}{4nt}>2q-\frac{(q-n-1)(q+2p-2)}{n}\label{9}.
\end{align}
For each $j\in\{1, \cdots, q\}$ and $k\in\{p,\cdots,t\},$ set $A_j^k:=\{z: \nu_{(f,H_j)}(z)=k\}$ and $B_j^k:=\{z: \nu_{(g,H_j)}(z)=k\}.$ Then we have $\overline{ A_j^k}\backslash A_j^k\subseteq \text{sing} f^{-1}(H_j)$, where the closure is taken with respect to the usual topology and 
$ \text{sing} f^{-1}(H_j)$ means the singular locus of the (reduction of the) analytic set $f^{-1}(H_j)$ of codimension one. Indeed, otherwise there existed
$a\in\overline{ A_j^k}\backslash A_j^k\cap \text{reg} f^{-1}(H_j).$ Then $p_0:=\nu_{(f,H_j)}(a)\ne k.$ Since $a$ is a regular point of $f^{-1}(H_j),$  by the R\"{u}ckert Nullstellensatz (see \cite{CAS}) we can choose nonzero holomorphic functions $h,u$ on a neigborhood $U$ of $a$ such that $dh$ and $u$ have no zero point and $(f, H_j) = h^{p_0}.u$ on $U$. Since $a \in \overline{A_j^k}$, there exists $b \in A_{j}^k \cap U$. Then $k = \nu_{(f, H_j)}(b) = \nu_{h^{p_0}.u}(b) = p_0$. This is a contradiction. Thus, $\overline{A_{j}^k} \backslash A_{j}^k \subseteq \text{sing}f^{-1}(H_j),$ for all $j \in  \{1, \cdots, q\}$ and $ k \in \{p,\dots, t\}$. This means that $\overline{A_{j}^k} \backslash A_{j}^k$ is a closed subset of an analytic set of codimension $\geq 2$. On the other hand $A_{j}^k \cap A_{j}^l =\varnothing$ for all $p\leq k \ne l\leq t$. Hence
\begin{align*}
(n-p)N^{[1]}(r, \overline{A_{j}^p}) +  \dots +(n-t)N^{[1]}(r, \overline{A_{j}^t}) \leq nN_{(f, H_j)}^{[1]}- N_{(f, H_j)}^{[n]}(r)\text{\ and\ }\\
p(n-p)N^{[1]}(r, \overline{A_{j}^p}) +  \dots +p(n-t)N^{[1]}(r, \overline{A_{j}^t}) \leq nN_{(f, H_j)}^{[p]}- pN_{(f, H_j)}^{[n]}(r)
\end{align*}
 for all  $j\in \{1, \cdots, q\}$ (note that $p\leq t< n).$
This implies that
\begin{align*}
\frac{n-t}{n}\sum_{k=p}^{t}N^{[1]}(r,\overline{A_{j}^k})\leq N_{(f, H_j)}^{[1]}(r)- \frac{1}{n}N_{(f, H_j)}^{[n]}(r)\text{\ and\ }\\
\frac{p(n-t)}{n}\sum_{k=p}^{t}N^{[1]}(r,\overline{A_{j}^k})\leq N_{(f, H_j)}^{[p]}(r)- \frac{p}{n}N_{(f, H_j)}^{[n]}(r),
\end{align*}
for all  $j\in \{1, \cdots, q\}.$
This gives
\begin{align}
\frac{(n-t)(q+2p-2)}{n}\sum_{j=1}^{q}\sum_{k=p}^{t}N^{[1]}(r,\overline{A_{j}^k})\quad\quad\quad\quad\quad\quad\quad\quad\notag\\ 
\leq (q-2)\sum_{j=1}^{q}\big(N_{(f, H_j)}^{[1]}(r)
- \frac{1}{n}N_{(f, H_j)}^{[n]}(r)\big)
+2\sum_{j=1}^{q}\big(N_{(f, H_j)}^{[p]}(r)- \frac{p}{n}N_{(f, H_j)}^{[n]}(r)\big)\label{10}.
\end{align}
Similarly,
\begin{align}
\frac{(n-t)(q+2p-2)}{n}\sum_{j=1}^{q}\sum_{k=p}^{t}N^{[1]}(r,\overline{B_{j}^k})\quad\quad\quad\quad\quad\quad\quad\quad\notag\\ 
\leq (q-2)\sum_{j=1}^{q}\big(N_{(g, H_j)}^{[1]}(r)
- \frac{1}{n}N_{(g, H_j)}^{[n]}(r)\big)
+2\sum_{j=1}^{q}\big(N_{(g, H_j)}^{[p]}(r)- \frac{p}{n}N_{(g, H_j)}^{[n]}(r)\big)\label{11}.
\end{align}
By (\ref{8}), (\ref{10}) and (\ref{11}) we have
\begin{align}
\frac{(n-t)(q+2p-2)}{n}\sum_{j=1}^{q}\sum_{k=p}^{t}\big( N^{[1]}(r,\overline{A_{j}^k})+N^{[1]}(r,\overline{B_{j}^k})\big)\quad\quad\quad\notag\\
\leq \big(2q-\frac{(q-n-1)(q+2p-2)}{n}\big)(T_f(r)+T_g(r))+o(T_f(r)+T_g(r))\label{12}.
\end{align}
Set $S_j^k:=A_{j}^k\cup B_{j}^k$ $(j\in\{1,\cdots,q\}, k\in\{p,\cdots,t\}).$

\noindent It is clear that
\begin{align*}
\min\{\nu_{(f,H_{j})},t+1\}=\min\{\nu_{(g,H_{j})},t+1\} \text{\ on\ } \C^m\backslash (\cup_{k=p}^t S_{j}^k)
\end{align*}
(note that $\min\{\nu_{(f,H_{j})},p\}=\min\{\nu_{(g,H_{j})},p\}$ on $\C^m$).

\noindent This means that
\begin{align*}
\{z: \min\{\nu_{(f,H_{j})}(z),t+1\}\ne\min\{\nu_{(g,H_{j})}(z),t+1\}\}\subset \cup_{k=p}^t S_{j}^k
\end{align*}
for all $j\in\{1,\cdots,q\}.$

\noindent Thus, for each $i\in\{1,\cdots,q\},$ by Lemma 6 (with $i_0=i, j_0=\sigma(i)$ and $\ell=t+1)$ we have
\begin{align*}
\sum_{j=1,j\ne i,\sigma(i)}^{q}N
_{(f,H_j)}^{[1]}(r)&+N_{(f,H_i)}^{[t+1]}(r)+ N_{(f,H_{\sigma(i)})}^{[t+1]}(r) \leq 
T_f(r)+T_g(r)\\
&+t\big(N^{[1]}(r,\overline{\cup_{k=p}^t S_{i}^k})
+N^{[1]}(r,\overline{\cup_{k=p}^t S_{\sigma(i)}^k})\big)+O(1).
\end{align*}
Then
\begin{align*}
&\sum_{i=1}^q\sum_{j=1,j\ne i,\sigma(i)}^{q}N
_{(f,H_j)}^{[1]}(r) +\sum_{i=1}^q\big(N_{(f,H_i)}^{[t+1]}(r)+ N_{(f,H_{\sigma(i)})}^{[t+1]}(r)\big)\\ 
&\leq 
q(T_f(r)+T_g(r))
+t\sum_{i=1}^q\big(N^{[1]}(r,\cup_{k=p}^t \overline{S_{i}^k})
+N^{[1]}(r,\cup_{k=p}^t \overline{S_{\sigma(i)}^k})\big)+O(1).
\end{align*}
On the other hand $\sigma:\{1,\cdots,q\}\to \{1,\cdots,q\}$ is bijective. Hence, we get
\begin{align}
(q-2)&\sum_{j=1}^qN
_{(f,H_j)}^{[1]}(r) +2\sum_{i=1}^qN_{(f,H_i)}^{[t+1]}(r)\notag\\
&\leq q(T_f(r)+T_g(r))
+2t\sum_{i=1}^qN^{[1]}(r,\cup_{k=p}^t \overline{S_{i}^k})+O(1)\notag\\
&\leq 
q(T_f(r)+T_g(r))
+2t\sum_{i=1}^q\sum_{k=p}^t\big (N^{[1]}(r, \overline{A_{i}^k}\cup\overline{B_{i}^k})\big)+O(1)\notag\\
&\leq 
q(T_f(r)+T_g(r))
+2t\sum_{i=1}^q\sum_{k=p}^t\big(N^{[1]}(r, \overline{A_{i}^k})+N^{[1]}(r,\overline{B_{i}^k})\big)+O(1).\label{13}
\end{align}
By the Second Main Theorem we have
\begin{align}
(q-2)\sum_{j=1}^qN
_{(f,H_j)}^{[1]}(r) &+2\sum_{i=1}^qN_{(f,H_i)}^{[t+1]}(r)\notag\\
&\geq \frac{(q-2)}{n}\sum_{j=1}^qN
_{(f,H_j)}^{[n]}(r) +\frac{2(t+1)}{n}\sum_{i=1}^qN_{(f,H_i)}^{[n]}(r)\notag\\
&\geq \frac{(q+2t)(q-n-1)}{n}T_f(r)-o(T_f(r))\label{14}
\end{align}
(note that $t+1\leq n).$

\noindent By (\ref{13}) and (\ref{14}) we have
\begin{align*}
\frac{(q+2t)(q-n-1)}{n}T_f(r)&-q\big(T_f(r)+T_g(r)\big)-o(T_f(r))\\
&\leq 2t\sum_{i=1}^q\sum_{k=p}^t\big(N^{[1]}(r, \overline{A_{i}^k})+N^{[1]}(r,\overline{B_{i}^k})\big).
\end{align*}
Similarly,
\begin{align*}
\frac{(q+2t)(q-n-1)}{n}T_g(r)&-q\big(T_f(r)+T_g(r)\big)-o(T_g(r))\\
&\leq 2t\sum_{i=1}^q\sum_{k=p}^t\big(N^{[1]}(r, \overline{A_{i}^k})+N^{[1]}(r,\overline{B_{i}^k})\big).
\end{align*}
Then
\begin{align*}
\big(\frac{(q+2t)(q-n-1)}{n}&-2q\big)(T_f(r)+T_g(r))-o(T_f(r)+T_g(r))\\
&\leq 4t\sum_{i=1}^q\sum_{k=p}^t\big(N^{[1]}(r, \overline{A_{i}^k})+N^{[1]}(r,\overline{B_{i}^k})\big).
\end{align*}
Combining with (\ref{12}) we get
\begin{align*}
&\big(\frac{(q+2t)(q-n-1)}{n}-2q\big)\frac{(n-t)(q+2p-2)}{4nt}(T_f(r)+T_g(r))\\
&\quad \leq \big(2q-\frac{(q-n-1)(q+2p-2)}{n}\big)(T_f(r)+T_g(r))+o(T_f(r)+T_g(r)).
\end{align*}
This implies that
\begin{align*}
\big(\frac{(q+2t)(q-n-1)}{n}-2q\big)\frac{(n-t)(q+2p-2)}{4nt}\leq 2q-\frac{(q-n-1)(q+2p-2)}{n}.
\end{align*}
This contradicts to (\ref{9}). Hence $f\equiv g.$ We have completed the proof of Theorem 3.\hfill$\Box$

 \noindent Gerd Dethloff \\
Universit\'{e} de Bretagne Occidentale \\
  UFR Sciences et Techniques \\
D\'{e}partement de Math\'{e}matiques \\
6, avenue Le Gorgeu, BP 452 \\
  29275 Brest Cedex, France \\
e-mail: gerd.dethloff@univ-brest.fr\\
\vskip0.15cm
      
 \noindent Tran Van Tan\\
Department of Mathematics\\
 Hanoi National University of 
Education\\
136 Xuan Thuy street, Cau Giay, Hanoi, Vietnam\\
e-mail:tranvantanhn@yahoo.com

\end{document}